\documentclass[11pt]{amsart}
\usepackage{amscd, amsmath, amsthm, amssymb}
\newcommand{\calC}{\mathcal{C}}

\newcommand{\calO}{\mathcal{O}}

\newcommand{\bbC}{\mathbb{C}}

\newcommand{\bbP}{\mathbb{P}}
\newcommand{\bbQ}{\mathbb{Q}}

\newcommand{\bbZ}{\mathbb{Z}}

\newcommand{\Aut}{\textup{Aut}}

\newcommand\Pic{{\text{Pic}}}

\newcommand\PU{{\textup{PU}}}

\DeclareMathOperator{\disc}{disc}

\newtheorem{theorem}{Theorem}[section]
\newtheorem{lemma}[theorem]{Lemma}
\newtheorem{proposition}[theorem]{Proposition}

\newtheorem{corollary}[theorem]{Corollary}
\theoremstyle{definition}     % italic or bold etc.
\newtheorem{definition}[theorem]{Definition}

\theoremstyle{remark}
\newtheorem{remark}[theorem]{Remark}
\numberwithin{equation}{section}
\begin{document}
\title[Quotients of fake projective planes]
{Quotients of fake projective planes}

\author[J. Keum]{JongHae Keum }
\address{School of Mathematics, Korea Institute for Advanced
Study, Seoul 130-722, Korea } \email{jhkeum@kias.re.kr}
\thanks{Research supported by the Korea Research Foundation Grant(KRF-2007-C00002).}

%\author[G. Prasad]{Gopal Prasad }
%\address{University of Michigan, Ann Arbor, MI 48109, USA } \email{gprasad@umich.edu}

\begin{abstract}
Recently, Prasad and Yeung classified all possible fundamental
groups of fake projective planes. According to their result, many
fake projective planes admit a nontrivial group of automorphisms,
and in that case it is isomorphic to $\bbZ/3\bbZ$, $\bbZ/7\bbZ$,
$7:3$, or $(\bbZ/3\bbZ)^2$, where $7:3$ is the unique non-abelian
group of order $21$.

Let $G$  be a group of automorphisms of a fake projective plane
$X$. In this paper we classify all possible structures of the
quotient surface $X/G$ and its minimal resolution.
\end{abstract}
\maketitle
%%%%%%%%%%%%%%%%%%%%%%%%%%%%%%%%%%%%%%%%%%%
%\setcounter{section}{0}
\section{Introduction}
It is known that a compact complex surface with the same Betti numbers as the
complex projective plane $\bbC\bbP^2$ is projective (see e.g. \cite{BHPV}).
Such a surface is called {\it a fake projective plane} if it is not
isomorphic to $\bbC\bbP^2$.

Mumford \cite{Mum} first proved the existence of a fake projective
 plane, based on the theory of the $p$-adic unit ball by Kurihara \cite{Ku} and
Mustafin \cite {Mus}. Later, using a similar idea, Ishida and Kato
\cite{IshidaKato} proved the existence of at least two more. Then,
Keum \cite{K} gave a construction of a fake projective plane with
an order 7 automorphism, using Ishida's description \cite{Ishida}
of an elliptic surface covered by a (blow-up) of Mumford's fake
projective plane. Recently, Prasad and Yeung \cite{PY} classified
all possible fundamental groups of fake projective planes.
According to their result, Keum's fake projective plane and
Mumford's fake projective plane are different from each other, but
belong to the same class. Furthermore, a group of automorphisms of
a fake projective plane is isomorphic to $\{1\}$, $\bbZ/3\bbZ$,
$\bbZ/7\bbZ$, $7:3$, or $(\bbZ/3\bbZ)^2$, and many fake projective
planes admit a nontrivial automorphism. Here
$7:3=(\bbZ/7\bbZ)\rtimes (\bbZ/3\bbZ)$ is the unique non-abelian
group of order $21$.

Let $G$ be a group of automorphisms of a fake projective plane
$X$.  In this paper we classify all possible structures of the
quotient surface $X/G$ and its minimal resolution. We first deal
with the case where $G$ is of prime order, and prove the
following:

\begin{theorem} \label{main}  Let $G$  be a group of automorphisms of a fake projective plane $X$. Let $Z= X/G$, and
$\nu : Y\to Z$ be a minimal resolution. Then the following two
statements are true.
\begin{enumerate}
\item[(1)] If the order of $G$ is $3$, then  $Z$ has $3$ singular points of type $\frac{1}{3}(1,2)$,
and  $Y$ is a minimal surface of general type with $K_{Y}^{2}=3, \,\,p_{g}=0$.
\item[(2)] If the order of $G$ is $7$, then $Z$ has $3$ singular points of type
$\frac{1}{7}(1,3)$, and $Y$ is a minimal elliptic surface of
Kodaira dimension $1$ with $2$ multiple fibres.  The pair of the
multiplicities is one of the following three cases: $(2, 3)$, $(2,
4)$, $(3, 3)$.
\end{enumerate}
\end{theorem}

We remark that the fundamental group $\pi_{1}(Y)$ of $Y$ is given
by $\{1\}$, $\bbZ/2\bbZ$, $\bbZ/3\bbZ$ in the three cases of (2),
respectively. (See \cite{D} for fundamental groups of elliptic
surfaces.) The first case of (2), where $Y$ is called a Dolgachev
surface, is supported by the example from \cite{K}. I have learnt
from Donald Cartwright and Tim Steger that according to their
computer calculation an order 7 quotient of a fake projective
plane has fundamental group either $\{1\}$ or $\bbZ/2\bbZ$. This
implies that the second case of (2) is supported by an example,
while the third case of (2) is not.

\begin{corollary} \label{maincor1} Let $X$ be a fake projective plane with $\Aut(X)\cong (\bbZ/3\bbZ)^2$.
Let $G=\Aut(X)$, $Z= X/G$, and $\nu : Y\to Z$ be a minimal
resolution. Then  $Z$ has $4$ singular points of type
$\frac{1}{3}(1,2)$, and $Y$ is a numerical Godeaux surface, i.e. a
minimal surface of general type with $K_{Y}^{2}=1, \,\,p_{g}=0$.
\end{corollary}

The above corollary proves the existence of a numerical Godeaux
surface with a configuration of 8 smooth rational curves of Dynkin
type $4A_2$, which has not been known before.

\begin{corollary} \label{maincor2} Let $X$ be a fake projective plane with $\Aut(X)\cong 7:3$.
Let $G=\Aut(X)$, $W= X/G$, and $\nu : V\to W$ be a minimal
resolution. Then  $W$ has $3$ singular points of type
$\frac{1}{3}(1,2)$ and $1$ singular point of type
$\frac{1}{7}(1,3)$. Furthermore, $V$ is a minimal elliptic surface
of Kodaira dimension $1$ with $2$ multiple fibres, and with $4$
reducible fibres of type $I_3$. The pair of the multiplicities is
the same as that of the minimal resolution of the order $7$
quotient of $X$.
\end{corollary}

\begin{corollary} \label{maincor3} Let $X$ be a fake projective plane with $\Aut(X)\cong 7:3$.
Let $G\cong \bbZ/7\bbZ<\Aut(X)$, $Z= X/G$, and $\nu : Y\to Z$ be a
minimal resolution. Then the elliptic fibration of $Y$ has $3$
singular fibres of type $I_1$, and $1$ reducible fibre of type
$I_9$.
\end{corollary}

Finally, we remark that there has not been known yet a geometric
construction of a fake projective plane that does not use the ball
quotient construction. A main purpose of this paper is to provide
some useful hints on how to find such a construction.

\bigskip
{\bf Acknowledgements.} I thank Gopal Prasad, Sai-Kee Yeung,
Donald Cartwright for many helpful conversations, and especially
Tim Steger for informing me of Lemma \ref{lem2}. I also thank Igor
Dolgachev for his comments, which helped me to improve the
exposition of the paper.

\bigskip
{\bf Notation}

\medskip\noindent
$D_{1} \equiv D_{2}$ : two divisors $D_{1}$ and $D_{2}$ are
linearly equivalent.\\
$D_{1} \sim D_{2}$ : two $\bbQ$-divisors $D_{1}$ and $D_{2}$ are
numerically equivalent.\\
$K_{X}$ : the canonical divisor of $X$.\\
$p_{g}(X)=\dim H^2(X, \calO_X)$ : the geometric genus of $X$.\\
$q(X)=\dim H^1(X, \calO_X)$ : the irregularity of $X$.\\
$\chi(X)=1-q(X)+p_q(X)$ : the holomorphic Euler characteristic.\\
$e(X)$ : the Euler number of $X$.\\
$b_i(X)$ : the i-th Betti number of $X$.\\
$c_i(X)$ : the i-th Chern class of $X$. $c_2(X)=e(X)$ if $X$ is a
smooth surface.\\
$g(C)$ : the genus of a curve $C$.\\
$(-n)$-curve : a smooth rational curve with self-intersection
$-n$.\\
singular point of type $\frac{1}{m}(1,a)$ : a cyclic quotient
singularity given by the diagonal action of diag$(\zeta,
\zeta^a)$ on $\bbC^2$, where $\zeta$ is a primitive $m$-th root of 1.\\
$\bbQ$-homology $\bbC\bbP^{2}$ : a normal projective surface with
the same Betti numbers as $\bbC\bbP^2$.
\medskip
\section{Preliminary Results}

There have been known many equivalent characterizations of a fake projective plane.

\begin{theorem} \label{def} A smooth compact complex surface $X$ with $b_1(X)=0$, $b_2(X)=1$ is a fake projective plane if one of the following holds true:
\begin{enumerate}
\item[(1)] $X$ is not isomorphic to $\bbC\bbP^2$.
\item[(2)] $X$ is not homeomorphic to $\bbC\bbP^2$.
\item[(3)] $X$ is not homotopy equivalent to $\bbC\bbP^2$.
\item[(4)]  $\pi_{1}(X)$ is an infinite group.
\item[(5)]  The universal cover of $X$ is a $2$-dimensional complex ball $B\subset \bbC^{2}$, and $X\cong B/\pi_{1}(X)$, where $\pi_{1}(X)\subset \PU(2,1)$.
\item[(6)] $K_{X}$ is ample.
\item[(7)] $K_{X}$ is ample, $p_{g}(X)=q(X)=0$, and $K_{X}^{2}=3c_{2}(X)=9$.
\end{enumerate}
\end{theorem}

We need the following lemma.

\begin{lemma} \label{curveonfpp}
Let $X$ be a fake projective plane, and $C$ be a smooth curve on $X$.
Then $e(C)\le -4$, or equivalently $g(C)\ge 3$.
\end{lemma}

\begin{proof}
Let $l$ be an ample generator of $\Pic(X)$ modulo torsions. Then $l^2=1$ and $C\equiv_{\bbQ} ml$ for some positive integer $m$.
Since $K_X\equiv_{\bbQ} 3l$, we have
$$e(C)=2-2g(C)=-C^2-CK_X= -(m^2+3m)\le -4.$$
\end{proof}

A normal projective complex surface is called a $\bbQ$-homology
$\bbC\bbP^{2}$ if it has the same Betti numbers with the complex
projective plane $\mathbb{C}\mathbb{P}^2$. If a $\bbQ$-homology
$\bbC\bbP^{2}$ is nonsingular, then it is either $\bbC\bbP^{2}$ or
a fake projective plane.

\begin{proposition} \label{rhp}
Let $S$ be a $\bbQ$-homology $\bbC\bbP^{2}$  with quotient
singularities only. Suppose that $S$ admits a finite group $G$ of
automorphisms. Then the quotient $S/G$ is again a $\bbQ$-homology
$\bbC\bbP^{2}$  with quotient singularities only.\\
In particular, $p_{g}(S/G)=q(S/G)=0$, $e(S/G)=3$ and
$\chi(S/G)=1$.
\end{proposition}

\begin{proof}
Since $S$ has $p_{g}=q=0$, so does the quotient $S/G$. Thus the
minimal resolution of $S/G$ has $q=0$, and hence $b_1=0$. It
follows that $b_1(S/G)=0$.

Since $S$ has $b_2=1$, so does the quotient $S/G$.
\end{proof}

 Now we consider fake
projective planes with an automorphism,  and get the following
preliminary information.

\begin{proposition} \label{quot}
Let $X$ be a fake projective plane with an automorphism $\sigma$.
Assume that the order of $\sigma$ is a prime number, say, $p$. Let
$Z=X/\langle\sigma\rangle$ and $\nu : Y\to Z$ be a minimal
resolution. Then
\begin{enumerate}
\item[(1)] $Z$ is a  $\bbQ$-homology $\bbC\bbP^{2}$ with $K_{Z}$ ample.
\item[(2)] $p_{g}(Y)=q(Y)=0$.
\item[(3)] $K_{Z}^{2}=\frac{9}{p}$.
\item[(4)]  The fixed point set $X^{\sigma}$ consists of $3$ points.
\end{enumerate}
\end{proposition}

\begin{proof}
The statements (1) and (2) follow from Proposition \ref{rhp}, as
$X$ itself is a  $\bbQ$-homology $\bbC\bbP^{2}$ with $K_{X}$
ample.

Note that $Z$ is a $\bbQ$-factorial variety, $K_{Z}$ is a
$\bbQ$-Cartier divisor, and
$$K_X\sim\pi^*K_{Z}$$
where $\pi:X\to Z$ is the quotient map.

Since $9=K_{X}^{2}=(\pi^*K_{Z})^2=pK_{Z}^2$, (3) follows. Here and
hereafter we use the intersection theory for $\bbQ$-divisors on
$\bbQ$-factorial varieties (or, for topologists, on V-manifolds).

It remains to prove (4). Assume that $X^{\sigma}$ contains a
curve. Since any two curves on a fake projective plane intersect,
we may assume that $X^{\sigma}$ consists of a smooth curve $C$ and
$r$ isolated points. Then the quotient surface $Z$ has $r$
singular points. Note that $e(X^{\sigma})=r+e(C)$. Using Hurwitz
formula, we have
$$e(X)=p\cdot e(Z) -(p-1)\{r+e(C)\}.$$
Since $e(X)=e(Z)=3$, this yields
\begin{equation}\label{1}
r+e(C)=3.
\end{equation}
This also follows from the topological Lefschetz fixed point
formula, as $\sigma$ acts identically on $H^*(X, \bbQ)$.

From the orbifold Bogomolov-Miyaoka-Yau inequality, one sees that
a $\bbQ$-homology $\bbC\bbP^{2}$ with quotient singularities only
cannot have more than 5 singular points (see e.g. \cite{Kollar},
\cite{HK}). Thus $r\le 5$. This bound together with Lemma
\ref{curveonfpp} contradicts \eqref{1}. Thus, $X^{\sigma}$ does
not contain a curve, and consists of $3$ isolated  points.
\end{proof}

\begin{corollary} \label{pneq2}
 $p\neq 2$.
\end{corollary}

\begin{proof}
Suppose $p=2$. Then $Z$ has rational double points only, hence
$K_{Y}\sim\nu^{*}K_{Z}$ and $K_{Y}^{2}=K_{Z}^{2}=\frac{9}{2}$,
which is not an integer, a contradiction.
\end{proof}

Prasad and Yeung \cite{PY} have provided precise possible values
for the order $p$. According to their result, $p=3$ or $7$. In
each case we will determine the types of singularities of the
quotient surface $Z$, using the holomorphic Lefschetz fixed point
formula.

\begin{lemma} \label{hlf} Let $S$ be a complex manifold of dimension $2$ with $p_g=q=0$.
Assume that $S$ admits an automorphism $\sigma$ of prime order
$p$. Let $r_i\,\, (1\le i\le p-1)$ be the number of isolated fixed
points of $\sigma$ which give singularities of type
$\frac{1}{p}(1,i)$ on the quotient surface. Let $C_1,\ldots, C_k$
be $1$-dimensional components of the fixed locus $S^{\sigma}$.
Then $$1=
\sum_{j=1}^k\Big\{\frac{1-g(C_j)}{2}+\frac{(p+1)C_j^2}{12}
\Big\}+\sum_{i=1}^{p-1}a_ir_i,$$ where
$$a_i=\frac{1}{p-1}\sum_{j=1}^{p-1}\frac{1}{(1-\zeta^j)(1-\zeta^{ij})}$$
with $\zeta$ a primitive pth root of $1$, e.g.
$a_1=\frac{5-p}{12},\,\,a_2=\frac{11-p}{24}$.
\end{lemma}

\begin{proof} This formula easily follows from the original holomorphic Lefschetz fixed point
formula \cite{ASIII}, p. 567. See, e.g. \cite{Zhang}, Lemma 1.6,
whose proof works for all complex manifolds of dimension $2$ with
$p_g=q=0$, not just for rational surfaces.
\end{proof}

\medskip
\section{The Case : $p=3$}

In this section we prove the following:

\begin{proposition} \label{p=3}  Let $\sigma$  be an automorphism of
 order $3$ of a fake projective plane $X$. Let
$Z=X/\langle\sigma\rangle$ and $\nu : Y\to Z$ be a minimal
resolution. Then $Z$ has $3$ singularities of type
$\frac{1}{3}(1,2)$, and $Y$ is a minimal surface of general type
with $K_{Y}^{2}=3, \,\,p_{g}=0$.
\end{proposition}

\begin{proof} From Proposition \ref{quot},
we see that $p_{g}(Y)=q(Y)=0$ and $K_{Z}^{2}=3$. Also we know that
the fixed locus $X^{\sigma}$ consists of $3$ points. Assume that
$Z$ has $r_i$ singular points of type $\frac{1}{3}(1,i)$. By Lemma
\ref{hlf} we have
$$1=a_1r_1+a_2r_2=\frac{1}{6}r_1+\frac{1}{3}r_2.$$
Since $r_1+r_2=3$, we see that $r_1=0$ and $r_2=3$. This proves
that $Z$ has $3$ singular points of type $\frac{1}{3}(1,2)$. In
particular, $K_{Y}\sim\nu^{*}K_{Z}$. Since $K_Z$ is ample, $K_Y$
is nef. Thus $Y$ is minimal.
 \end{proof}

\begin{corollary} \label{Z/3^2} Let $X$ be a fake projective plane with $\Aut(X)\cong (\bbZ/3\bbZ)^2$.
Let $G=\Aut(X)$, $Z= X/G$, and $\nu : Y\to Z$ be a minimal
resolution. Then  $Z$ has $4$ singular points of type
$\frac{1}{3}(1,2)$, and $Y$ is a minimal surface of general type
with $K_{Y}^{2}=1, \,\,p_{g}=0$.
\end{corollary}

\begin{proof} The group $G$ has 4 subgroups isomorphic to
$\bbZ/3\bbZ$. Each fixes 3 isolated points of type
$\frac{1}{3}(1,2)$. No stabilizer of a nonsingular point can be
isomorphic to $(\bbZ/3\bbZ)^2$, thus there are 12 points, each of
whose stabilizers is isomorphic to $\bbZ/3\bbZ$. It follows that
$Z$ has $4$ singular points of type $\frac{1}{3}(1,2)$. Note that
the canonical divisor $K_Z$ of $Z$ is $\bbQ$-Cartier and
$$K_X\sim\pi^*K_{Z}$$
where $\pi:X\to Z$ is the quotient map. Thus
$$K_Z^2=K_X^2/9=1.$$ Since $Z$ has only rational double points and $K_Z$ is ample, we
see that $K_Y$ is nef, $K_Y^2=1$ and hence the assertion on $Y$
follows.
\end{proof}

According to \cite{PY}, many fake projective planes admit an
automorphism of
 order $3$, thus by taking a quotient one can obtain many new examples of a minimal surface of general type with $K_{Y}^{2}=3, \,\,p_{g}=0$.

\bigskip\noindent
{\bf Question:} Does there exist a fake projective plane $X$ with
an automorphism $\sigma$ of
 order $3$ such that the minimal resolution $Y$ of $Z=X/\langle\sigma\rangle$ is simply
 connected?

\begin{remark} (1) Since $Z$ has rational singularities only, $\pi_1(Y)\cong \pi_1(Z)$.
Thus the question is whether there is a fake projective plane $X$
with an automorphism  $\sigma$  of
 order $3$ such that the augmented fundamental group $\langle\pi_{1}(X),\,\tilde{\sigma} \rangle$, where $\tilde{\sigma}$ is a lift of $\sigma$  to the ball, is
 the normal closure of the subgroup generated by elements of order 3.

 (2) In the case of order 7, there is a fake projective plane $X$ with an automorphism  $\sigma$  of
 order $7$ such that the minimal resolution $Y$ of the quotient is simply connected \cite{K}. In this case $Y$ is not of general type.

 (3) A simply connected surface of general type with $K^{2}=3, \,\,p_{g}=0$ has been recently constructed by H. Park, J. Park, and D. Shin \cite{PPS}.
 They use the method of Lee and Park \cite{LP}, which produces examples of simply connected surfaces of general type with $K^{2}=2, \,\,p_{g}=0$.
 An affirmative answer to the question would give yet another interesting example of
 a simply connected surface of general type with $K^{2}=3, \,\,p_{g}=0$.
 \end{remark}

\medskip
\section{The Case : $p=7$}

In this section we first prove the following:

\begin{proposition} \label{1,3only}  Let $\sigma$  be an automorphism of
 order $7$ of a fake projective plane $X$. Let
$Z=X/\langle\sigma\rangle$ and $\nu : Y\to Z$ be a minimal
resolution. Then $Z$ has $3$ singular points of type
$\frac{1}{7}(1,3)$, and $K_{Y}^{2}=0$.
\end{proposition}

This follows from the following three lemmas.

\begin{lemma} \label{disect}  Let $\sigma$  be an automorphism of
 order $7$ of a fake projective plane $X$. Let
$Z=X/\langle\sigma\rangle$ and $\nu : Y\to Z$ be a minimal
resolution. Then $Z$ has either $3$ singular points of type
$\frac{1}{7}(1,3)$, or  $2$ singular points of type
$\frac{1}{7}(1,4)$,
 and $1$ singular point of type $\frac{1}{7}(1,6)$.
\end{lemma}

\begin{proof} From Proposition \ref{quot},
we know that the fixed locus $X^{\sigma}$ consists of $3$ points.
Assume that $Z$ has $r_i$ singular points of type
$\frac{1}{7}(1,i)$. By Lemma \ref{hlf} we have
$$-r_1+r_2+2r_3+r_4+2r_5+4r_6=6.$$
Adding this to $\sum r_i=3$, we get
$$2(r_2+r_4)+3(r_3+r_5)+5r_6=9.$$
If $r_6=0$, then $r_3+r_5=3$, hence we get 3 points of type
$\frac{1}{7}(1,3)=\frac{1}{7}(1,5)$.\\ If $r_6=1$, then
$r_2+r_4=2$, hence we get 1 point of type $\frac{1}{7}(1,6)$ and
$2$ points of type $\frac{1}{7}(1,2)=\frac{1}{7}(1,4)$. This
proves the assertion.
\end{proof}

\begin{lemma} \label{lem2}  Let $\sigma$  be an automorphism of
 order $7$ of a fake projective plane $X$.
Then $\sigma$ cannot have a fixed point of type
$\frac{1}{7}(1,4)$.
\end{lemma}

\begin{proof} Assume that $\sigma$ has a fixed point of type
$\frac{1}{7}(1,4)$. Then the group acting on the complex ball
$B\subset \bbC^{2}$ must contain a matrix $\in \PU(2,1)$ which
diagonalises as:
\begin{displaymath}
M = \left( \begin{array}{ccc}
\alpha & 0 &0 \\
0 & \alpha\zeta & 0 \\
0 & 0 & \alpha\zeta^4 \\
\end{array} \right)
\end{displaymath}
where $\zeta=\zeta_7=e^{2\pi i/7}$ is the 7-th root of unity, and
$\alpha$ a complex number.

Using the notation of \cite{PY}, we can choose this matrix to be
in $\bar{\Gamma}$, which is contained in a rank 3 division algebra
over the field denoted by $\ell$. Therefore
  $${\rm tr}(M)= \alpha(1 + \zeta +
  \zeta^4)\,\,\, {\rm and}\,\,\, \det(M)= \alpha^3\zeta^5$$ both must belong to $\ell$. Thus
$\ell$ contains ${\rm tr}(M)^3/\det(M)$, which is equal to
$$(1+\zeta+\zeta^4)^3 / \zeta^5=6(\zeta+\zeta^{-1})^3+(\zeta+\zeta^{-1})^2-15(\zeta+\zeta^{-1})+5.$$
The field which this generates over $\bbQ$, namely
$\bbQ[\zeta+\zeta^{-1}]$, must be contained in $\ell$.  None of
the cases on Prasad-Yeung's final list has such an $\ell$. There
is exactly one possibility listed, but later excluded, which does
have such an~$\ell$, namely $\calC_{31}$.
\end{proof}

\begin{definition} \label{9curv} When $Z$ has $3$ singularities of type
$\frac{1}{7}(1,3)$, we denote by $$A_1, A_2, A_3, B_1, B_2, B_3,
C_1, C_2, C_3$$ the exceptional curves of $\nu:Y\to Z$ whose
Dynkin diagram is given by
$$(-2) \textrm{---}(-2)\textrm{---}(-3)\quad (-2)\textrm{---}(-2)\textrm{---}(-3)\quad (-2)
\textrm{---}(-2)\textrm{---}(-3)$$ \end{definition}

\begin{lemma} \label{7-1lem1}  Assume that $Z$ has $3$ singularities of type $\frac{1}{7}(1,3)$.
Then
$$K_Y\sim\nu^{*}K_{Z}-\frac{1}{7}(A_1+2A_2+3A_3)-\frac{1}{7}(B_1+2B_2+3B_3)-\frac{1}{7}(C_1+2C_2+3C_3).$$
In particular, $K_Y^2=0$.
\end{lemma}
\begin{proof} This follows from the
adjunction formula
$$K_{Y}\sim\nu^{*}K_{Z}-D,$$
where $D$ is a $\bbQ$-linear combination of the exceptional curves
$A_i,\,B_i,\, C_i$ with coefficients in the interval $[0, \,1)$.
These coefficients can be uniquely determined by the system of 9
linear equations
$$-2 -E^2=K_Y\cdot E=-D\cdot E\quad{\rm for}\quad E=A_i,\,B_i,\, C_i (i=1,2,3).$$
Finally note that $K_{Z}^{2}=\frac{9}{7}$ by Proposition
\ref{quot}. Also note that
$$\{\frac{1}{7}(A_1+2A_2+3A_3)\}^2=-K_Y\cdot\{\frac{1}{7}(A_1+2A_2+3A_3)\}=-\frac{3}{7}.$$ Hence the last
assertion follows.
\end{proof}

This completes the proof of Proposition \ref{1,3only}.\\
To complete the proof of Theorem \ref{main}, it suffices to prove
the following:

\begin{proposition} \label{p=7-1} Assume that $Z$ has $3$ singular points of
type $\frac{1}{7}(1,3)$. Then there are three cases :
\begin{enumerate}
\item[(i)]  $Y$ is a minimal elliptic surface of Kodaira dimension $1$ with $2$ multiple  fibres with multiplicity $2$ and $3$, respectively.
\item[(ii)]  $Y$ is a minimal elliptic surface of Kodaira dimension $1$ with $2$ multiple  fibres with multiplicity $2$ and $4$, respectively.
\item[(iii)]  $Y$ is a minimal elliptic surface of Kodaira dimension $1$ with $2$ multiple  fibres with multiplicity $3$ and $3$, respectively.
\end{enumerate}
\end{proposition}

The proof of Proposition \ref{p=7-1} consists of several lemmas.

\begin{lemma} \label{7-1lem2}  Assume that $Z$ has $3$ singularities of type $\frac{1}{7}(1,3)$.
Then
\begin{enumerate}
\item[(1)]  $-mK_Y$ is not effective for any positive integer $m$.
\item[(2)]  The Kodaira dimension of $Y$ is at least $1$.
\end{enumerate}
\end{lemma}

\begin{proof}
(1) Note that for $m\ge 1$
$$(\nu^{*}K_{Z})\cdot (-mK_Y)=-m(\nu^{*}K_{Z})^2=\frac{-9m}{7} < 0.$$
Since $\nu^{*}K_{Z}$ is nef, $-mK_Y$ cannot be effective.

(2) By Proposition \ref{quot}, $p_{g}(Y)=q(Y)=0$. By Proposition
\ref{1,3only}, $K_{Y}^{2}=0$. Thus if $Y$ is of Kodaira dimension
$\le 0$, then by the classification theory of surfaces $Y$ is
either a rational surface or an Enriques surface.

From the Riemann-Roch theorem and (1) we see that for any integer $k \ge 2$
$$h^0(kK_Y)=1+h^1(kK_Y)\ge 1.$$
Thus $Y$ is not rational.

Since $(\nu^{*}K_{Z})\cdot
(K_Y)=(\nu^{*}K_{Z})^2=K_{Z}^2=\frac{9}{7} > 0$, $K_Y$ is not
numerically trivial, hence $Y$ is not an Enriques surface. This
proves (2).
\end{proof}

Note that $\Pic(Y)\cong H^2(Y,\bbZ)$. Let $\Pic(Y)_f:=\Pic(Y)/{\rm
torsion}$. With the intersection pairing $\Pic(Y)_f$ becomes a
lattice.
\begin{lemma} \label{7-1lem3}  Assume that $Z$ has $3$ singularities of type
$\frac{1}{7}(1,3)$.  Then one can choose two $\bbQ$-divisors
$$\renewcommand{\arraystretch}{1.3}\begin{array}{l}
L=\frac{1}{7}(A_1+2A_2+3A_3)+\frac{2}{7}(B_1+2B_2+3B_3)+\frac{4}{7}(C_1+2C_2+3C_3),\\
M=\frac{1}{3}\nu^{*}K_{Z}-\frac{2}{7}(B_1+2B_2+3B_3)+\frac{1}{7}(C_1+2C_2+3C_3),
\end{array}$$
such that the lattice $\Pic(Y)_f$ is generated over the integers
by the numerical equivalence classes of $M$, $L$, and the $8$
curves $A_2, A_3, B_1, B_2, B_3, C_1, C_2, C_3.$
\end{lemma}

\begin{proof}
Note first that in this case $K_Y^2=0$, hence by Noether formula
rank$\Pic(Y)_f=10$. Since $\Pic(Y)_f$ contains an element of
self-intersection $-3$, e.g. $A_3$,  it is odd unimodular and of
signature $(1,9)$.

Let $R$ be the sublattice of $\Pic(Y)_f$ generated by the
numerical equivalence classes of the $9$ curves $A_1, A_2, A_3,
B_1, B_2, B_3, C_1, C_2, C_3.$ Let $\overline{R}$ and $R^{\perp}$
be its primitive closure and its orthogonal complement,
respectively, in the lattice $\Pic(Y)_f$. Note that $R^{\perp}$ is
of rank 1.\\ For an integral lattice $N$, let $\disc(N)$ denote
the discriminant group of $N$
$$\disc(N):=Hom(N, \bbZ)/N.$$
It is easy to see that $\disc(R)\cong (\bbZ/7\bbZ)^3$, more
precisely
$$\disc(R)=\langle\frac{1}{7}(A_1+2A_2+3A_3), \,\frac{1}{7}(B_1+2B_2+3B_3), \, \frac{1}{7}(C_1+2C_2+3C_3)\rangle.$$
Note that the length (the minimum number of generators) of
$\disc(R)$ is 3. Since the lattice $\Pic(Y)_f$ is unimodular,
$\disc(\overline{R})$ is isomorphic to $\disc(R^{\perp})$ which is
of length 1. Hence $R$ must be of index 7 in $\overline{R}$, and
the generator of $\overline{R}/R$ is of the form
$$L=\frac{1}{7}(A_1+2A_2+3A_3)+\frac{a}{7}(B_1+2B_2+3B_3)+\frac{b}{7}(C_1+2C_2+3C_3).$$
Since both $L\cdot K_Y$ and $L^2$ must be integers, we see that
$(a,b)=(2,4)$ or $(4,2)$ modulo 7. Thus up to interchanging the
curves $B_i$'s and $C_i$'s, we have determined the divisor $L$
uniquely modulo $R$.

Now we have  $\disc(\overline{R})\cong \disc(R^{\perp})\cong
\bbZ/7\bbZ$. Note that the integral divisor $7\nu^{*}K_{Z}$
belongs to $R^{\perp}$ and $(7\nu^{*}K_{Z})^2=7\cdot 3^2$. Thus
$R^{\perp}$ is generated by $\frac{7}{3}\nu^{*}K_{Z}$, hence
$$\disc(R^{\perp})=  \langle \frac{1}{3}\nu^{*}K_{Z}\rangle.$$
On the other hand,
$$\disc(\overline{R})=\langle L\rangle^{\perp}/\langle L\rangle=\langle \frac{3}{7}(B_1+2B_2+3B_3)+\frac{2}{7}(C_1+2C_2+3C_3)\rangle,$$
where $\langle L\rangle=\overline{R}/R$ is the isotropic subgroup
of $\disc(R)$ generated by $L($mod $R)$ and $\langle
L\rangle^{\perp}$ is its orthogonal complement in  $\disc(R)$ with
respect to the discriminant quadratic form on  $\disc(R)$. (See,
e.g. \cite{Nik} for discriminant quadratic forms for integral
lattices.) Thus the index 7 extension $\overline{R}\oplus
R^{\perp}\subset \Pic(Y)_f$ is given by the element of the form
$$M=\frac{1}{3}\nu^{*}K_{Z}+a\{\frac{3}{7}(B_1+2B_2+3B_3)+\frac{2}{7}(C_1+2C_2+3C_3)\}.$$
Since $M\cdot K_Y$ is an integer, we see that $a=4$ modulo 7.
This determines the divisor $M$ uniquely modulo $R$.
\end{proof}

\begin{remark}\label{intmat} The proof of Lemma
\ref{7-1lem3} shows that the two $\bbQ$-divisors $M$ and $L$ are
indeed integral divisors modulo torsion. The intersection matrix
of the 10-divisors $M, L, A_2, A_3, B_1, B_2, B_3, C_1, C_2, C_3$
is given by
\begin{displaymath} \left( \begin{array}{cccccccccc}-2& 0 &0 &0& 0& 0& 2& 0& 0& -1\\
 0& -9& 0& -1& 0& 0& -2& 0& 0& -4\\
0& 0 &-2 &1& 0& 0& 0& 0& 0& 0\\
 0& -1& 1& -3& 0& 0& 0& 0& 0& 0\\
 0& 0 &0 &0& -2& 1& 0& 0& 0& 0\\
 0& 0 &0 &0& 1& -2& 1& 0& 0& 0\\
 2& -2& 0& 0& 0& 1& -3& 0& 0& 0\\
 0& 0 &0 &0& 0& 0& 0& -2& 1& 0\\
 0& 0& 0& 0& 0& 0& 0& 1& -2& 1\\
 -1& -4 &0 &0& 0& 0& 0& 0& 1& -3\\
\end{array} \right)
\end{displaymath}
It is easy to see that this matrix has determinant $-1$. This
double checks that our choice of $M$, $L$ and 8 curves was
correct. But the unimodularity of the matrix only is not enough to
prove Lemma \ref{7-1lem3}, as it does not imply the integrality
(modulo torsion) of $M$ and $L$. There are many possible choices
of non-integral $\bbQ$-divisors $M$ and $L$ making the matrix
unimodular.
\end{remark}

\begin{lemma} \label{7-1lem4}  Assume that $Z$ has $3$ singularities of type $\frac{1}{7}(1,3)$.
Then $Y$ does not contain a $(-1)$-curve $($a smooth rational
curve of self-intersection $-1)$ $E$ with $0<E\cdot \nu^* K_Z<
\frac{9}{7}$.
\end{lemma}

\begin{proof} Assume that $Y$ contains such a $(-1)$-curve $E$.
Write
$$E\sim mM-dL+a_2 A_2+a_3 A_3+b_1 B_1+b_2 B_2+b_3 B_3+c_1 C_1+c_2 C_2+c_3 C_3$$
with integer coefficients. Since $E\cdot \nu^* K_Z=mM\cdot \nu^*
K_Z=\frac{3m}{7}$, the condition $0<E\cdot \nu^*K_Z<\frac{9}{7}$
is equivalent to $1\le m\le 2$.

\medskip
(1) Assume that $m=1$, i.e.
$$E\sim M-dL+a_2 A_2+a_3 A_3+b_1 B_1+b_2 B_2+b_3 B_3+c_1 C_1+c_2 C_2+c_3 C_3.$$
Then, the coefficients in the above satisfy the following system
of 9 inequalities and one equality:

$$\renewcommand{\arraystretch}{1.3}\begin{array}{l}
0\le E\cdot A_1=a_2\\
0\le E\cdot A_2=-2a_2+a_3\\
0\le E\cdot A_3=d+a_2-3a_3\\
0\le E\cdot B_1=-2b_1+b_2\\
0\le E\cdot B_2=b_1-2b_2+b_3\\
0\le E\cdot B_3=2+2d+b_2-3b_3\\
0\le E\cdot C_1=-2c_1+c_2\\
0\le E\cdot C_2=c_1-2c_2+c_3\\
0\le E\cdot C_3=-1+4d+c_2-3c_3\\
-1= E\cdot K_Y=-3d+a_3+b_3+c_3
\end{array}$$

From the 9 inequalities of the system, we obtain that
\begin{equation}\label{8}
a_3\le \frac{2}{5}d,\qquad b_3\le \frac{3}{7}(2+2d),\qquad c_3\le \frac{3}{7}(-1+4d).
\end{equation}
Indeed, eliminating $a_2$ from the second and the third inequality
of the system, we get the first inequality of \eqref{8}.
Eliminating $b_1$ and $b_2$ from the fourth, the fifth and the
sixth inequality of the system, we get the second inequality of
\eqref{8}. The third inequality of \eqref{8} can be proved
similarly.

Also, we obtain the following bound for $d$.
\begin{equation}\label{9}
0\le d\le 50.
\end{equation}
Indeed, from the first three inequalities of the system, we have
$$d\ge -a_2+3a_3=3(-2a_2+a_3)+5a_2\ge 5a_2\ge 0.$$
Applying the three inequalities of \eqref{8} to the equality of
the system, we get
$$3d-1=a_3+b_3+c_3\le \frac{2}{5}d + \frac{3}{7}(2+2d)+\frac{3}{7}(-1+4d),$$
hence $d\le 50$.

We know that $E^2=-1$. Expanding $E^2$ using the intersection
matrix from Remark \ref{intmat} and then applying the equality of
the system, we get
$$1+3d^2+2d= (4+2d)b_3+(6d-2)c_3+(a_2 A_2+a_3 A_3)^2+(\sum_{i=1}^3b_i B_i)^2+(\sum_{i=1}^3c_i C_i)^2.$$
Note that
$$\renewcommand{\arraystretch}{1.3}\begin{array}{l}
(a_2 A_2+a_3
A_3)^2=-2a_2^2+2a_2a_3-3a_3^2=-2(a_2-\frac{1}{2}a_3)^2-\frac{5}{2}a_3^2\le
-\frac{5}{2}a_3^2,\\
(\sum_{i=1}^3b_i
B_i)^2=-2(b_1-\frac{1}{2}b_2)^2-\frac{3}{2}(b_2-\frac{2}{3}b_3)^2-\frac{7}{3}b_3^2\le
-\frac{7}{3}b_3^2,\\
$$(\sum_{i=1}^3c_i C_i)^2=-2(c_1-\frac{1}{2}c_2)^2-\frac{3}{2}(c_2-\frac{2}{3}c_3)^2-\frac{7}{3}c_3^2\le -\frac{7}{3}c_3^2,
\end{array}$$
forcing the above equality to give the following inequality:
\begin{equation}\label{10}
1+3d^2+2d\le -\frac{5}{2}a_3^2-\frac{7}{3}b_3^2-\frac{7}{3}c_3^2+(4+2d)b_3+(6d-2)c_3.
\end{equation}

We claim that there is no integer solution satisfying
\eqref{8}-\eqref{10} and the equality of the system, hence no
solution satisfying the system.

The proof is cumbersome and goes as follows. First we obtain the
list of solutions $(d, a_3, b_3, c_3)$ of the equality of the
system under the constraints given by \eqref{9} and \eqref{8}; for
each value of $d$ from \eqref{9}, we solve the equation
$3d-1=a_3+b_3+c_3$ in the range \eqref{8}.\\ The following list is
generated by a computer program.
$$\renewcommand{\arraystretch}{1.3}\begin{array}{l}
(40,16,35,68), (33,13,29,56), (30,12,26,51), (26,10,23,44), (25,10,22,42),\\
(23,9,20,39), (20,8,18,33), (19,7,17,32), (18,7,16,30), (16,6,14,27),\\
(15,6,13,25), (13,5,12,21), (12,4,11,20), (11,4,10,18), (10,4,9,16),\\
(9,3,8,15), (8,3,7,13), (6,2,6,9), (5,2,5,7), (5,2,4,8),
(5,1,5,8),\\ (4,1,4,6), (3,1,3,4), (2,0,2,3), (1,0,1,1),
(0,0,0,-1)
\end{array}$$
Next, it is easy to check that none of these satisfies \eqref{10}.

\medskip
(2) Assume that $m=2$, i.e.
$$E\sim 2M-dL+a_2 A_2+a_3 A_3+b_1 B_1+b_2 B_2+b_3 B_3+c_1 C_1+c_2 C_2+c_3 C_3.$$
In this case the coefficients satisfy the following system of 9
inequalities and one equality:

$$\renewcommand{\arraystretch}{1.3}\begin{array}{l}
0\le E\cdot A_1=a_2\\
0\le E\cdot A_2=-2a_2+a_3\\
0\le E\cdot A_3=d+a_2-3a_3\\
0\le E\cdot B_1=-2b_1+b_2\\
0\le E\cdot B_2=b_1-2b_2+b_3\\
0\le E\cdot B_3=4+2d+b_2-3b_3\\
0\le E\cdot C_1=-2c_1+c_2\\
0\le E\cdot C_2=c_1-2c_2+c_3\\
0\le E\cdot C_3=-2+4d+c_2-3c_3\\
-1= E\cdot K_Y=-3d+a_3+b_3+c_3
\end{array}$$
Also in this case, \eqref{8}-\eqref{10} are replaced by
\begin{equation}\label{11}
a_3\le \frac{2}{5}d,\qquad b_3\le \frac{3}{7}(4+2d),\qquad c_3\le
\frac{3}{7}(-2+4d).
\end{equation}
\begin{equation}\label{12}
0\le d\le 65.
\end{equation}
\begin{equation}\label{13}
7+3d^2+2d\le
-\frac{5}{2}a_3^2-\frac{7}{3}b_3^2-\frac{7}{3}c_3^2+(8+2d)b_3+(6d-4)c_3.
\end{equation}

The same argument as in the case $m=1$ shows that there is no
solution satisfying \eqref{11}-\eqref{13} and the equality of the
system, hence no solution satisfying the system.
\end{proof}

\begin{lemma} \label{7-1lem5}  Assume that $Z$ has $3$ singularities of type $\frac{1}{7}(1,3)$.
Then $Y$ is minimal.
\end{lemma}

\begin{proof}
From Proposition \ref{1,3only} we know that $K_Y^2=0$.

Suppose $Y$ is not minimal. Then by Lemma \ref{7-1lem2} (2), $Y$
is of general type. Let $\mu : Y\to Y'$ be a birational morphism
to the minimal model $Y'$. Then
$$K_Y\equiv \mu^* K_{Y'} +\sum E_i,$$
where $E_i$'s are effective divisors, not necessarily irreducible,
with $E_i^2=-1$, $E_i\cdot E_j=0$ for $i\neq j$. Note that a
positive multiple of $\mu^* K_{Y'}\equiv K_Y-\sum E_i$ is
effective. Since $\nu^* K_Z$ is nef, we have $$(K_Y-\sum E_i)\cdot
\nu^* K_Z \ge 0.$$ Furthermore, $\mu^* K_{Y'}$ has positive
self-intersection, thus by Hodge index theorem $$(K_Y-\sum
E_i)\cdot \nu^* K_Z \neq 0.$$ Summarizing these, we have
\begin{equation}\label{18}
(K_Y-\sum E_i)\cdot \nu^* K_Z >0.
\end{equation}
Let $E$ be a $(-1)$-curve on $Y$. Since $\nu^* K_Z$ is nef and $E$
is not contracted by $\nu$, we have $$E\cdot \nu^* K_Z >0.$$ On
the other hand, by \eqref{18} we have $$E\cdot \nu^* K_Z<K_Y\cdot
\nu^* K_Z=\frac{9}{7}.$$ Thus, the assertion follows from Lemma
\ref{7-1lem4}.
\end{proof}

{\it Proof of Proposition \ref{p=7-1}.}

By Lemma \ref{7-1lem2} and \ref{7-1lem5}, $Y$ is a minimal
elliptic surface of Kodaira dimension $1$. It remains to prove the
assertion on multiplicities of multiple fibres.

Let $|F|$ be the elliptic pencil on $Y$. By the canonical bundle
formula for elliptic fibrations (see e.g. \cite{BHPV} Chap V),
$$F\sim nK_Y$$
for some positive rational number $n$. We claim that $n$ must be
an integer. To see this, we first note that $Y$ contains a
$(-3)$-curve, e.g. the curve $A_3$ (Definition \ref{9curv}), hence
$A_3\cdot K_Y=1$. Thus $n=A_3\cdot F$ is an integer.

Let $m_1F_1, m_2F_2, \dots, m_rF_r$ be the multiple fibres of the
elliptic fibration with multiplicity $m_1, m_2, \dots, m_r$,
respectively. Since $Y$ is not rational, $r\ge 2$. Again by the
canonical bundle formula for elliptic fibrations,
$$K_Y\equiv -F+\sum_{i=1}^{r}(m_i-1)F_i\equiv (r-1)F-\sum_{i=1}^{r}F_i,$$
hence
\begin{equation}\label{19}
\frac{1}{n}=r-1-\sum_{i=1}^{r}\frac{1}{m_i}.
\end{equation}
Since $\sum_{i=1}^{r}\frac{1}{m_i}\le \frac{r}{2}$, \eqref{19}
implies that  $r\le 3$ if $n=2$ ; $r=2$ if $n\ge 3$. Since
$n=A_3\cdot F=m_iA_3\cdot F_i$, all $m_i$ divide $n$. With these,
a further analysis of \eqref{19} shows that if $n=2$, then $(m_1,
m_2, m_3)=(2, 2, 2)$ ; if $n=3$, then $(m_1, m_2)=(3, 3)$ ; if
$n=4$, then $(m_1, m_2)=(2, 4)$ ; if $n=6$, then $(m_1, m_2)=(2,
3)$ ; if $n=5$ or $n\ge 7$, then there is no solution for $m_i$'s.
The first case would imply that there is a degree 2 map $:A_3\to
{\bbP}^1$, ramified at 3 points, which is impossible. This
completes the proof of Proposition \ref{p=7-1}.

\begin{corollary} \label{7:3} Let $X$ be a fake projective plane with $\Aut(X)\cong 7:3$.
Let $G=\Aut(X)$, $W= X/G$, and $\nu : V\to W$ be a minimal
resolution. Then  $W$ has $3$ singular points of type
$\frac{1}{3}(1,2)$ and $1$ singular point of type
$\frac{1}{7}(1,3)$. Furthermore, $V$ is a minimal elliptic surface
of Kodaira dimension $1$ with $2$ multiple fibres, and with $4$
reducible fibres of type $I_3$. The pair of the multiplicities is
the same as that of the minimal resolution of the order $7$
quotient of $X$.
\end{corollary}

\begin{proof} Write
$$G=\langle\sigma, \tau\,|\,\, \sigma^7=\tau^3=1,
\tau\sigma\tau^{-1}=\sigma^2\rangle.$$ Let
$Z=X/\langle\sigma\rangle$ be the order $7$ quotient of $X$, and
$Y$ be a minimal resolution of $Z$. Then by Proposition
\ref{p=7-1}, $Z$ has Kodaira dimension 1, and has $3$ singular
points of type $\frac{1}{7}(1,3)$, which form a single orbit of
the induced automorphism $\bar{\tau}$. By Proposition \ref{p=3},
every element of order 3 of $G$ fixes 3 points of type
$\frac{1}{3}(1,2)$. No stabilizer of a nonsingular point can be
isomorphic to $7:3$, thus $W=Z/\langle\bar{\tau}\rangle$ has $3$
singular points of type $\frac{1}{3}(1,2)$ and $1$ singular point
of type $\frac{1}{7}(1,3)$.\\ Note that the canonical divisor
$K_W$ of $W$ is an ample $\bbQ$-Cartier divisor and
$$K_W^2=\frac{K_X^2}{21}=\frac{3}{7}.$$ Thus by adjunction formula $K_V^2=0$. Since $Y$
has Kodaira dimension 1, $V$ has Kodaira dimension $\le 1$. Note
that the action of $\bar{\tau}$ on $Z$ lifts to $Y$. Let
$W'=Y/\langle\bar{\tau}\rangle$. By Proposition \ref{p=7-1}, we
know that $K_Y$ is nef. Thus $K_{W'}$ is nef. Since $W'$ has 3
singular points of type $\frac{1}{3}(1,2)$, and since $V$ is the
minimal resolution of $W'$, we see that $K_V$ is nef. This proves
that $V$ is minimal and of Kodaira dimension $\ge 0$. Note that
$$(\nu^{*}K_{W})\cdot (K_V)=K_{W}^2=\frac{3}{7}
> 0,$$ thus $K_V$ is not numerically trivial. This proves that $V$ has Kodaira dimension $1$.
The elliptic fibration on $V$ is given by a multiple of $K_V$.\\
Now $V$ has 9 smooth rational curves coming from the resolution
$\nu : V\to W$. The eight $(-2)$-curves among them must be
contained in fibres of the elliptic fibration. This is possible
only if the
fibres are the union of 4 reducible fibres of type $I_3$, since $V$ has Picard number 10.\\
Note that $Y$ is the degree 3 cover of $W'$ branched along the 3
singular points of $W'$, and $W'$ has an elliptic fibration
structure with a $(-3)$-curve that is a multi-section. The
$(-3)$-curve on $W'$ splits in $Y$ giving three $(-3)$-curves,
thus the elliptic fibres of $W'$ do not split in $Y$. The fibre
containing one of the singular point of $W'$ gives a fibre of type
$I_1$, the fibre of type $I_3$ gives a fibre of type $I_9$, and
the multiple fibres give multiple fibres of the same
multiplicities.
\end{proof}

\begin{remark} (1) By Proposition \ref{p=7-1}, the elliptic fibration on $V$ has multiplicities (2,3) or
(2,4) or (3,3). Such an elliptic surface with multiplicities (2,3)
was constructed by Ishida \cite{Ishida}. His construction was
based on the description of Mumford surface as a ball quotient.
From his elliptic surface a fake projective plane was constructed
\cite{K}.

(2) By Cartwright and Steger such an elliptic surface with
multiplicities (2,4) exists. From such an elliptic surface one can
give a similar construction of a fake projective plane \cite{K2}.
 \end{remark}

From the proof of Corollary \ref{7:3}, we also have the following:

\begin{corollary} \label{7:3'} Let $X$ be a fake projective plane with $\Aut(X)\cong 7:3$.
Let $G\cong \bbZ/7\bbZ<\Aut(X)$, $Z= X/G$, and $\nu : Y\to Z$ be a
minimal resolution. Then the elliptic fibration of $Y$ has $3$
singular fibres of type $I_1$, and $1$ reducible fibre of type
$I_9$.
\end{corollary}

%%%%%%%%%%%%%% reference %%%%%%%%%%%%%%%%%


\begin{thebibliography}{[BPV]}
\bibitem[ASIII]{ASIII} M. F. Atiyah, and I. M. Singer,
\textit{The index of elliptic operators, III}, Ann. of Math. {\bf
87} (1968), 546-604.
\bibitem[BHPV]{BHPV} W. Barth, K. Hulek, Ch. Peters, and A. Van de Ven,
\textit{Compact Complex Surfaces}, second ed. Springer 2004.
\bibitem[D]{D} I. Dolgachev, \textit{Algebraic surfaces with
$q=p_g=0$}, in CIME "Algebraic Surfaces", Liguori Edit., Naples
(1981), 97-216.
\bibitem[HK]{HK} D. Hwang, and J. Keum, \textit{The maximum number of singular points on rational homology projective planes},
math.AG/0801.3021.
\bibitem[Is]{Ishida}  M. Ishida, \textit {An elliptic surface covered
    by Mumford's fake projective plane}, Tohoku Math. J. {\bf 40}
(1988), 367-398.
\bibitem[IsKa]{IshidaKato}  M. Ishida, F. Kato, \textit {The strong
    rigidity
theorem for non-archimedean uniformization}, Tohoku Math. J. {\bf
50} (1998), 537-555.
\bibitem[K]{K} J. Keum, \textit{A fake projective plane with an order 7 automorphism}, Topology {\bf 45} (2006),
919-927.
\bibitem[K2]{K2} J. Keum, \textit{A fake projective plane constructed from an elliptic surface with multiplicities $(2,4)$}, preprint.
\bibitem[Ko]{Kollar} J. Koll\'ar, \textit{Is there a topological Bogomolv-Miyaoka-Yau inequality?},\\ math.AG/0602562.
%\bibitem[KZ]{KZ} J. Keum, and D.-Q. Zhang, \textit{Algebraic surfaces with quotient singularities -
%including some discussion on automorphisms and fundamental groups}, Proceedings of
%``Algebraic Geometry in East Asia'' Kyoto, 2001,
\bibitem[Ku]{Ku} A. Kurihara, \textit {Construction of $p$-adic
    unit balls and Hirzebruch proportionality}, Amer. J. Math. {\bf 102}
(1980), 565-648.
\bibitem[LP]{LP} Y. Lee, and J. Park, \textit{ A simply connected surface of general type with $p_{g}=0$ and $K^{2}=2$}, Invent. Math.
{\bf 170} (2007), 483-505.
\bibitem[Mum]{Mum} D. Mumford, \textit{An algebraic surface with $K$
    ample $K^2=9, p_g=q=0$}, Amer. J. Math. {\bf 101} (1979),
233-244.
\bibitem[Mus]{Mus} G.A. Mustafin, \textit {Non-archimedean
    Uniformization},
Math. USSR, Sbornik, {\bf 34}
(1978), 187-214.
\bibitem[N]{Nik} V. V. Nikulin, \textit{Integral symmetric bilinear forms and its applications}, Izv. Akad. Nauk SSSR Ser. Mat. \textbf{43} (1979), no. 1,
111-177; English translation: Math. USSR Izv. \textbf{14} (1979),
no. 1, 103-167 (1980)

\bibitem[PPS]{PPS} H. Park, J. Park, and D. Shin, \textit{ A simply connected surface of general type with $p_{g}=0$ and $K^{2}=3$}, math.AG/07080273.
\bibitem[PY]{PY} G. Prasad, and S.-K. Yeung, \textit{Fake projective planes}, Invent. Math. {\bf 168} (2007), 321-370.
\bibitem[Z]{Zhang} D.-Q. Zhang, \textit{Automorphisms of finite order on rational surfaces}, J. Algebra {\bf 238} (2001), 560-589.
\end{thebibliography}
\end{document}